\documentclass[11pt]{article}
\usepackage{amssymb,amsmath,amsthm}
\usepackage{subfigure}

\title{A note on packing spanning trees in graphs\\ and bases in matroids}
\author{Robert~F.~Bailey\footnote{Division of Science (Mathematics), Grenfell Campus, Memorial University of \mbox{Newfoundland,} University Drive, Corner Brook, Newfoundland A2H~6P9, Canada. E-mail: \texttt{rbailey@grenfell.mun.ca}},\hspace{2mm}
Mike~Newman\footnote{Department of Mathematics and Statistics, University of Ottawa, 585 King Edward Avenue, Ottawa, \mbox{Ontario} K1N~6N5, Canada. E-mail: \texttt{mnewman@uottawa.ca}}\hspace{2mm} and
Brett~Stevens\footnote{School of Mathematics and Statistics, Carleton University, 1125 Colonel By Drive, Ottawa, Ontario K1S~5B6, Canada.  E-mail: \texttt{brett@math.carleton.ca}}}

\newtheorem{thm}{Theorem}
\newtheorem{lemma}[thm]{Lemma}
\newtheorem{prop}[thm]{Proposition}

\theoremstyle{definition}
\newtheorem{defn}[thm]{Definition}
\newtheorem{example}[thm]{Example}

\newtheorem*{notation}{Notation}

\newcommand{\blob}{\circle*{0.2}}
\newcommand{\complete}{
  \begin{picture}(2,2)
  \put(0,0){\blob}
  \put(2,0){\blob}
  \put(0,2){\blob}
  \put(2,2){\blob}

  \put(0,0){\line(1,0){2}}
  \put(0,0){\line(0,1){2}}
  \put(0,2){\line(1,0){2}}
  \put(2,0){\line(0,1){2}}
  \put(0,0){\line(1,1){2}}
  \put(0,2){\line(1,-1){2}}
  \end{picture} }

\DeclareMathOperator{\rank}{r}

\begin{document}

\maketitle

\begin{abstract}
We consider the class of graphs for which the edge connectivity is equal to the maximum number of edge-disjoint spanning trees, and the natural generalization to matroids, where the cogirth is equal to the number of disjoint bases.  We provide descriptions of such graphs and matroids, showing that such a graph (or matroid) has a unique decomposition.  In the case of graphs, our results are relevant for certain communication protocols.
\end{abstract}



\section{Introduction} \label{section:intro}

In \cite{Itai88}, Itai and Rodeh proposed a communication protocol, called the {\em $k$-tree protocol}, which allows all nodes of a network to communicate through a distinguished root node~$v$, even when some set of $k-1$ or fewer edges are removed from the network.  The protocol requires the graph $G$ modelling the network to have two properties.  First, the graph $G$ must remain connected when any $k-1$ edges are removed, so $k$ can be at most the edge connectivity of $G$.  Second, it requires a collection of $k$ spanning trees for $G$, $\{T_1,\ldots,T_k\}$, with the following property (which they called the {\em $k$-tree condition for edges}): for all vertices $w$ distinct from $v$, and for any $i,j$ where $1\leq i < j \leq k$, the paths in $T_i$ and $T_j$ from $v$ to $w$ are internally disjoint.
\newpage

Clearly, if $G$ has $k$ edge-disjoint spanning trees, then it satisfies the $k$-tree condition.  (This is not a requirement: for example, a cycle satisfies the $2$-tree condition, but does not have two edge-disjoint spanning trees.)  In particular, if the number of edge-disjoint spanning trees (denoted $\sigma(G)$) is equal to the edge connectivity (denoted $\lambda(G)$), then Itai and Rodeh's protocol can be applied; we call such graphs {\em maximum spanning tree packable}, or {\em max-STP}.

In~\cite{ubbnet}, two of the present authors considered a related network protocol, where they require a collection $\mathcal{U}$ of spanning trees (not necessarily pairwise disjoint) for $G$ chosen so that for any $t$ edges (where $t<\lambda(G)$), there exists a spanning tree $T \in \mathcal{U}$ disjoint from those $t$ edges.  Ideally, this collection of (not necessarily disjoint) spanning trees (called an {\em uncovering-by-bases}, or UBB, for $G$) should be as small as possible.  The class of max-STP graphs with $\lambda(G)=\sigma(G)=k$ is also of interest here, as the $k$ edge-disjoint spanning trees form a UBB which is optimal in two ways: (i) the spanning trees are disjoint, so the UBB is as small as possible; and (ii) the number of edges which can be ``uncovered'' by the collection is as large as possible.  (In fact, it was the study of UBBs which led the authors to the results in the present paper.)  The notion of UBBs also generalizes to matroids: see \cite[Section~7]{btubb}.  Another class of matroids where UBBs arise naturally are as follows.

A {\em base} for a group acting on a set is a subset of points whose pointwise stabilizer is trivial; equivalently, every group element is uniquely specified by its action on those points.  (See~\cite{Cameron99} for more details.)  A UBB for a permutation group is a collection of bases so that any $t$-subset of points is disjoint from some base in the collection; these have applications to the decoding of permutation codes (see~\cite{ecpg}).  In~\cite{CameronFDF95}, Cameron and Fon-Der-Flaass investigated permutation groups whose bases form the bases of a matroid: such groups are known as {\em IBIS groups}.  In the case of an IBIS group, a UBB for the group is also a UBB for the corresponding matroid.  An important sub-class of IBIS groups are the {\em base-transitive} groups, where the bases lie in a single orbit of the group; constructions of UBBs for many examples of base-transitive groups are given in~\cite{btubb}.

A straightforward example of a max-STP graph $G$ with $\lambda(G)=\sigma(G)=k$ is the graph obtained from a tree $T$ with $n$ edges by replacing each edge of $T$ by $k$ parallel edges.  A UBB for $G$ can be formed from disjoint copies of $T$.  The cycle matroid of $G$ is obtained from a free matroid by replacing each point with a parallel class of size $k$ (regardless of the structure of $T$); alternatively, this is the transversal matroid of a uniform set-partition.  This same matroid arises from a base-transitive group: if $H$ is a group acting regularly on a set $X$ of $k$ points, consider the wreath product $H\wr S_n$ acting on $n$ disjoint copies of $X$, labelled $X_1,\ldots,X_n$.  A base for this group consists of a single point chosen from each of $X_1,\ldots,X_n$, and the corresponding matroid is again the transversal matroid of the set system $\{X_1,\ldots,X_n\}$.  As a code, this group can correct $d=\lfloor (k-1)/2 \rfloor$ errors, and a UBB requires $d+1$ disjoint transversals; in the cycle matroid of $G$, this is equivalent to $d+1$ disjoint spanning trees.  (This class of codes is discussed in more detail in~\cite{thomas}.)

\subsection{Graphs}
Let $G=(V,E)$ be a graph.  Throughout this paper, we shall assume that graphs are connected, and we allow for the possibility of multiple edges.

An {\em edge cut} in $G$ is a partition $(V_1,V_2)$ of the vertex set of $G$ into two non-empty subsets.  The corresponding set of edges are those with one endpoint in $V_1$ and the other endpoint in $V_2$.  Removing the edges of an edge cut disconnects $G$; if the number of such edges is $k$, we call it a {\em $k$-edge cut}.  In a mild abuse of terminology, we will use the term ``edge cut'' to refer both to the partition of the vertex set and the set of edges of the cut.  The {\em edge-connectivity} of $G$, denoted $\lambda(G)$, is the least value of $k$ for which there exists a $k$-edge cut in $G$.  We note that sometimes we will refer to a $k$-edge cut by the set of edges whose removal disconnects the graph, rather than the partition of $V$.  Also, we say that $G$ is {\em $k$-edge connected} if $\lambda(G) \geq k$.

The {\em spanning tree packing number} of $G$, denoted $\sigma(G)$, is the maximum number of edge-disjoint spanning trees in $G$.  (We usually shorten this to {\em STP number}.)  A survey of results on STP numbers can be found in Palmer \cite{Palmer01}.  In particular, graphs with given STP number were characterized independently by Nash-Williams \cite{NashWilliams61} and Tutte \cite{Tutte61}, both in 1961.  

\begin{thm}[Nash-Williams; Tutte] \label{thm:TutteNW}
A connected graph $G$ has at least $k$ edge-disjoint spanning trees if and only if, for every partition of $V(G)$ into $r$ parts, there are at least $k(r-1)$ edges between the parts.
\end{thm}

It is a straightforward observation that $\sigma(G) \leq \lambda(G)$: clearly, to disconnect $G$ we must remove at least one edge from each of the $\sigma(G)$ disjoint spanning trees (and possibly some other edges as well).  Also, in 1983 Gusfield \cite{Gusfield83} showed that it follows from Nash-Williams and Tutte's result that $\lambda(G) \leq 2\sigma(G)$ (see also Diestel \cite[Section 3.5]{Diestel00}).  When presented with an inequality such as $\sigma(G) \leq \lambda(G)$, it seems natural to ask when equality is achieved, $\sigma(G) =  \lambda(G)$.

\begin{defn} \label{defn:maxSTP}
A graph $G$ is said to be {\em maximum spanning tree-packable}, or {\em max-STP} for short, if $\lambda(G)=\sigma(G)$, i.e.\ the edge connectivity is equal to the spanning tree packing number.
\end{defn}

\begin{example} \label{example:easy}
The graph~$G$ in Figure~\ref{fig:maxSTPexample} is a max-STP graph with $\lambda(G)=\sigma(G)=2$.  (In fact, as $G$ has $8$ vertices and $14$ edges, any pair of edge-disjoint spanning trees must contain all edges of $G$.)
\setlength{\unitlength}{7mm}
\begin{figure}[h]
\centering
\begin{picture}(7,2.2)
\put(-0.15,0){\complete}

\put(4.85,0){\complete}

\put(2,0){\line(1,0){3}}
\put(2,2){\line(1,0){3}}
\end{picture}
\caption{A max-STP graph with $\sigma(G)=\lambda(G)=2$. \label{fig:maxSTPexample}}
\end{figure}
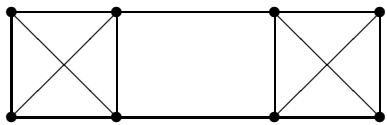
\end{example}

In Section~\ref{section:graphdecomp}, we will present a structure theorem for max-STP graphs.  The two parameters $\lambda(G)$ and $\sigma(G)$ both have straightforward analogues in matroid theory, so in Section~\ref{section:matroiddecomp} we prove an analogous theorem for matroids where the two parameters agree.

\subsection{Matroids}

A {\em matroid} $M$ consists of a (finite) ground set $E$ together with a family $\mathcal{I}$ of subsets of $E$, called {\em independent sets}, which satisfy the following three axioms:
\begin{itemize}
\item[I1.] $\mathcal{I} \neq \emptyset$;
\item[I2.] if $I\in \mathcal{I}$ and $J\subseteq I$, then $J\in \mathcal{I}$;
\item[I3.] if $I,J\in \mathcal{I}$ and $|J|<|I|$, then there exists $x \in I\setminus J$ such that $J\cup\{x\} \in \mathcal{I}$.
\end{itemize}

The maximal independent sets are called the {\em bases} of $M$; the collection of these is denoted $\mathcal{B}$.  The bases are necessarily equicardinal; this size is called the {\em rank} of the matroid and is denoted by $\rank(M)$.  The rank of an arbitrary subset $X \subseteq E$ is the size of the largest independent set contained in $X$ and is denoted by $\rank(X)$.  For any subsets $X$ and $Y$, the rank satisfies the inequality $\rank(X \cap Y) \leq \rank(X)+\rank(Y)-\rank(X \cup Y)$. A {\em flat} of rank $k$ in a matroid is a maximal set of rank $k$.  The intersection of flats is always a flat; a {\em hyperplane} is a flat of rank $\rank(M)-1$.  A {\em cocircuit} is a minimal set that intersects every basis, i.e.\ a minimal subset $S \subseteq E$ for which $S\cap B_i \neq \emptyset$ for all $B_i \in \mathcal{B}$.  Equivalently, the cocircuits of a matroid are exactly the complements of the hyperplanes of the matroid.  If $M=(E,\mathcal{I})$ is a matroid and $X\subseteq E$, the matroid obtained by the {\em deletion} of $X$, denoted $M\setminus X$, has ground set $E\setminus X$, and its collection of independent sets is $\{ I\setminus X \, : \, I\in\mathcal{I} \}$.  (For background material on matroids, we refer the reader to Oxley \cite[Chapter~1]{Oxley92}).

One of the motivating examples of matroids (and the one most relevant to this paper) is the {\em cycle matroid} of a graph $G=(V,E)$, where the ground set $E$ is indeed the edge set of $G$, and where the bases are the maximum spanning forests of $G$.  We denote this matroid by $M(G)$.  The independent sets of $M(G)$ are the subsets of $E$ which contain no cycles of $G$.  In particular if $G$ is a connected graph with $n$ vertices, the bases of $M(G)$ are the spanning trees of $G$, and $M(G)$ has rank $n-1$.  More generally, if $G$ has $c$ connected components, then $M(G)$ has rank $n-c$.  A cocircuit of $M(G)$ is a minimal set of edges whose removal increases the number of components of $G$ by one, i.e., a minimal edge-cut.

The natural analogue of spanning tree packing number for matroids is as follows.

\begin{defn} \label{defn:bpn}
Let $M=(E,\mathcal{I})$ be a matroid.  The {\em base packing number} of $M$, denoted $\sigma(M)$, is the size of the largest set of disjoint bases of $M$.
\end{defn}

Matroids with given base packing number were described by Edmonds \cite{Edmonds65} in 1965 (see also Oxley \cite{Oxley92}, Theorem 12.3.11), thereby generalizing the result of Nash-Williams and Tutte (Theorem~\ref{thm:TutteNW} above).  Edmonds showed that a matroid $M=(E,\mathcal{I})$ has $k$ disjoint bases if and only if, for every subset $X \subseteq E$, the following inequality holds:
\[ k\cdot\rank(X) + |E\setminus X| \geq k\cdot\rank(M).\]

Edge-connectivity also has a natural analogue for matroids.

\begin{defn} \label{defn:cogirth}
Let $M=(E,\mathcal{I})$ be a matroid.  The {\em cogirth} of $M$, denoted $\lambda(M)$, is the smallest size of a cocircuit in $M$.
\end{defn}

The inequality $\sigma(G)\leq\lambda(G)$ carries over to matroids in a straightforward way.

\begin{prop} \label{prop:ineq}
Let $M=(E,\mathcal{I})$ be a matroid.  Then $\sigma(M)\leq \lambda(M)$.
\end{prop}

\proof Given any collection of $\sigma(M)$ disjoint bases, any cocircuit of $M$ must intersect each of these bases in at least one element. \endproof

As with graphs, when presented with the inequality in Proposition \ref{prop:ineq}, it seems natural to ask when equality is achieved.  In Section~\ref{section:matroiddecomp}, we present a result which describes the structure of matroids for which $\sigma(M)=\lambda(M)$.

\subsection{Preliminary results} \label{subsection:prelims}

We begin with a brief discussion of matroid connectivity.  A matroid $M=(E,\mathcal{I})$ is {\em disconnected} if there is a partition of $E$ into non-empty sets $E_1,E_2$, such that there are matroids $M_1=(E_1,\mathcal{B}_1)$ and $M_2=(E_2,\mathcal{B}_2)$ with the property that every base of $M$ is the union of a base for $M_1$ and a base for $M_2$.  In this situation, we write $M= M_1\oplus M_2$.  If there is no such partition, we say $M$ is {\em connected}.

Now, for the cycle matroid of a graph $G$, it is easy to see that if $G$ is a disconnected graph, then $M(G)$ is disconnected as a matroid.  However, the converse is not true in general: if $G$ is connected but contains a cut vertex, $M(G)$ is disconnected.  The appropriate matroid definition of being connected in the graphic matroid is in fact equivalent to the graph being 2-connected. The fact that the notions of connectivity do not coincide is the main reason why we present separate analyses for graphs and matroids in this paper.  In particular, Theorem~\ref{thm:main} for graphs is not simply a corollary of Theorem~\ref{thm:decomp} for matroids.  (For a further discussion of matroid connectivity, see Oxley~\cite[Chapter 4]{Oxley92}, and \cite[\S8.2]{Oxley92} for a comparison of this notion with that of graph connectivity.)

The following lemma is straightforward, but will prove to be crucial to us.

\begin{lemma} \label{lemma:disjointcocircuits}
Suppose that $M=(E,\mathcal{I})$ is a matroid for which $\sigma(M)=\lambda(M)=k$.  Then any pair of distinct minimum cocircuits of $M$ are disjoint.
\end{lemma}

\proof Suppose for a contradiction that $C_1$ and $C_2$ are distinct minimum cocircuits of $M$ contain some element $x$ in common.  Now, $E\setminus (C_1 \cup C_2) = (E\setminus C_1) \cap (E \setminus C_2)$ is the intersection of two distinct hyperplanes, and therefore has rank at most $\rank(M)-2$ ($\ast$).

Consider a collection of $k=\sigma(M)$ disjoint bases $B_1,\ldots,B_k$; note that $x\in B_i$ for some $i$, since each cocircuit intersects each basis non-trivially.  Since $|C_1|=|C_2|=\lambda(M)=\sigma(M)=k$, it follows that each of the $k$ bases contains exactly one element of $C_1$ and of $C_2$.  In particular, $B_i \cap (C_1 \cup C_2) = \{x\}$, so $B_i \setminus \{x\} \subseteq E\setminus (C_1 \cup C_2)$.  But the rank of $B_i \setminus \{x\}$ is $\rank(M)-1$, contradicting ($\ast$).

Hence $C_1 \cap C_2 = \emptyset$.  \endproof

In the case of graphs, Lemma~\ref{lemma:disjointcocircuits} states that in a max-STP graph, no pair of minimal $k$-edge cuts can have an edge in common.  This can be shown directly by a straightforward counting argument, similar to that of Zhang~\cite[Lemma 2.2]{Zhang02}.

The next lemma is also not difficult.

\begin{lemma} \label{lemma:deletecocircuit}
Suppose that $\sigma(M)=k$.  Then for any cocircuit $C$ of $M$ of size $k$, and where $C\neq E$, we have $\sigma(M\setminus C)\geq k$.
\end{lemma}

\proof Let $\{B_1,\ldots,B_k\}$ be $k$ disjoint bases of $M$.  Each $B_i$ intersects $C$ in a unique element, say $x_i$.  For $1 \leq i \leq k$, let $B_i'=B_i\setminus \{x_i\}$.  Then it follows that $\{B_1',\ldots,B_k'\}$ is a set of disjoint bases for $M\setminus C$.
\endproof

In the case of graphs, Lemma~\ref{lemma:deletecocircuit} states that in a graph with $k$ edge-disjoint spanning trees, if an $k$-edge cut $C$ is deleted, then both connected components of the resulting graph $G\setminus C$ will also have at least $k$ edge-disjoint spanning trees (unless that component is an isolated vertex).

From now on, we will consider the cases of graphs and matroids separately.


\section{Decomposing graphs} \label{section:graphdecomp}
In this section, we will obtain our structural description of the max-STP graphs.  To assist with this, we define the following ``joining'' operation.

\begin{defn} \label{defn:kjoin}
Let $G_1=(V_1,E_1)$ and $G_2=(V_2,E_2)$ be connected graphs, where $V_1 \cap V_2 =\emptyset$, and let $K$ be a set of $k$ edges with one end in $V_1$ and one end in $V_2$ (for some integer $k$).  Then the {\em $K$-join} of $G_1$ and $G_2$, denoted by $G_1 \ast_K G_2$, is the graph with vertex set $V_1 \cup V_2$ and edge set $E_1 \cup E_2 \cup K$.
\end{defn}

When the set $K$ of $k$ edges is not specified (or is not important), we speak of a {\em $k$-join} of $G_1$ and $G_2$, and denote it by $G_1 \ast_k G_2$.

We follow the definition with a couple of remarks.  First, by construction $(V_1,V_2)$ is a $k$-edge cut of $G_1 \ast_K G_2$, so consequently the edge-connectivity of the $K$-join is at most $k$.  Second, two $k$-joins will not, in general, be isomorphic (unless we have a special case, such as when $k=1$ and both $G_1$ and $G_2$ are vertex-transitive).  Example~\ref{example:K4*2K4} below shows the kind of situation which may arise.

\begin{example} \label{example:K4*2K4}
There are three non-isomorphic possibilities for $K_4 \ast_2 K_4$ (the $2$-join of two copies of $K_4$), as shown in Figure \ref{fig:K4*2K4}.
\end{example}
\setlength{\unitlength}{7mm}
\begin{figure}[hbtp]
\centering
\begin{picture}(7,8)
\put(-0.15,6){\complete}

\put(4.85,6){\complete}

\put(2,6){\line(1,0){3}}
\put(2,8){\line(1,0){3}}

\put(-0.15,3){\complete}

\put(4.85,3){\complete}

\put(2,5){\line(1,0){3}}
\put(2,5){\line(3,-2){3}}

\put(-0.15,0){\complete}

\put(4.85,0){\complete}

\qbezier(2,2)(3.5,2.5)(5,2)
\qbezier(2,2)(3.5,1.5)(5,2)

\end{picture}
\caption{The three non-isomorphic possibilities for $K_4 \ast_2 K_4$. \label{fig:K4*2K4}}
\end{figure}
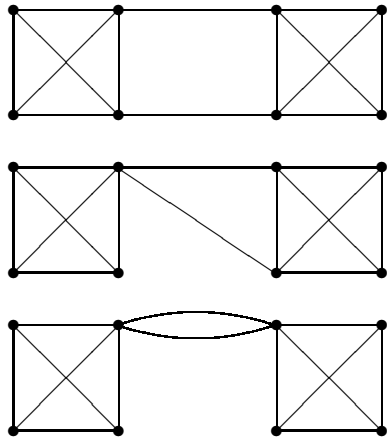

We are interested in $k$-joins because they preserve the property we are concerned with.  First, it is a straightforward exercise to show that if a graph $G$ has the form $G=G_1 \ast_K G_2$, where $k = |K| \leq \sigma(G_i) \leq \lambda(G_i)$ (for $i=1,2$), then $\lambda(G)=\sigma(G)=k$.  Conversely, if $G$ is a max-STP graph with $\lambda(G)=\sigma(G)=k$, then $G$ is necessarily a $k$-join of graphs $G_1$, $G_2$, and where $G_1$, $G_2$ each satisfy exactly one of the following:
\begin{itemize}
  \item[(i)] $G_i$ has one vertex and no edges;
  \item[(ii)] $k \leq \sigma(G_i) < \lambda(G_i)$;
  \item[(iii)] $k < \sigma(G_i) = \lambda(G_i)$;
  \item[(iv)] $\sigma(G_i) = \lambda(G_i) = k$.
\end{itemize}

\begin{defn} \label{defn:kirred}
We call a graph {\em $k$-irreducible} if it belongs to classes (i)--(iii) above; we call it {\em $k$-reducible} if it belongs to class (iv).
\end{defn}

We remark that these four classes (i)--(iv) partition the class of all graphs with $k$ or more edge-disjoint spanning trees (i.e.~the class determined by Nash-Williams and Tutte).  Also, we observe that a graph $G$ in class (iii) will itself be a max-STP graph, but with a higher spanning tree packing number and edge-connectivity; such a graph will also be $k'$-reducible, where $k'=\sigma(G)=\lambda(G)>k$.

A complication arises when two or more $k$-joins are made.  Suppose that we make a $k$-join $H=G_1 \ast_k G_2$, and then make $H \ast_k G_3$ (where $G_i=(V_i,E_i)$).  If the $k$ edges in the second $k$-join are all attached to exactly one of $G_1$ or $G_2$ (assume without loss of generality that this is $G_2$), then in the resulting graph $(V_1,V_2\cup V_3)$ and $(V_1\cup V_2,V_3)$ will both be $k$-edge cuts, obtained by removing the edges of the first or second $k$-joins respectively.  However, if the second $k$-join attaches $G_3$ to some vertices in each of $G_1$ and $G_2$, only $(V_1\cup V_2,V_3)$ is a $k$-edge cut.  This phenomenon is demonstrated in Example \ref{example:orderindep} below.

\begin{example} \label{example:orderindep}
Figure \ref{figure:orderindep} shows two ways of forming 2-joins of three copies of $K_4$.  The graph on the left has two 2-edge cuts, while the graph on the right has only one.
\end{example}
\setlength{\unitlength}{7mm}
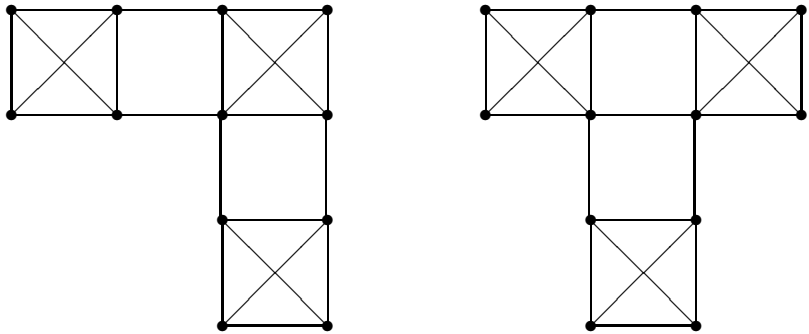
\begin{figure}[hbtp]
\centering
\begin{picture}(15,6)
\put(-0.15,4){\complete}

\put(3.85,4){\complete}

\put(2,4){\line(1,0){2}}
\put(2,6){\line(1,0){2}}

\put(3.85,0){\complete}

\put(4,2){\line(0,1){2}}
\put(6,2){\line(0,1){2}}


\put(8.85,4){\complete}

\put(12.85,4){\complete}

\put(11,4){\line(1,0){2}}
\put(11,6){\line(1,0){2}}

\put(10.85,0){\complete}

\put(11,2){\line(0,1){2}}
\put(13,2){\line(0,1){2}}

\end{picture}
\caption{Two ways of forming $2$-joins of three copies of $K_4$. \label{figure:orderindep} }
\end{figure}

We call a sequence of $k$-joins {\em order-independent} if the joining edges of any one of them yield a $k$-edge cut in the resulting graph.  Thus in Figure \ref{figure:orderindep}, the 2-joins in the graph on the left are order-independent, while those in the graph on the right are not.

The term order-independent refers to the fact that the final graph is invariant of different choices of which $k$-join is performed first, second, third, etc.  However, which ``pieces'' are joined by a particular $k$-join remain fixed.

In the case where $G$ is a max-STP graph with $\lambda(G)=\sigma(G)=k$, and which has more than one $k$-edge cut, Lemma~\ref{lemma:disjointcocircuits} tells us that no pair of $k$-edge cuts can ``overlap'' (i.e.\ they can have no edge in common).  Furthermore, by construction, if we have made two $k$-joins order-independently, then the two $k$-edge cuts arising from these must be non-overlapping.  Following this, we are now able to state our decomposition theorem for max-STP graphs.

\begin{thm} \label{thm:main}
Suppose that $G$ is a max-STP graph satisfying $\lambda(G)=\sigma(G)=k$.  Then we have the following.
\begin{itemize}
\item[(i)] There exists a unique set $\mathcal{A}$ of $k$-irreducible graphs $G_1,\ldots,G_m$ (for some $m$).
\item[(ii)] There exists a unique rooted tree $R$ with $m$ leaves labelled by $G_1,\ldots,G_m$, such that the root is labelled by $G$ and each node is labelled by an order-independent $k$-join of its children.
\item[(iii)] For each non-leaf, labelled by $H$ and its $d$ children labelled $H_1,\ldots,H_d$, there exists a unique tree $T_H$ with vertices $\{1,\ldots,d\}$ labelled by $H_1,\ldots,H_d$, such that for each edge $e=ij$ of $T_H$, there exists a $k$-edge cut $K_e$ of $H$ such that $H_i \ast_{K_e} H_j$ is an induced subgraph of $H$.
\end{itemize}
\end{thm}

We remark that Theorem~\ref{thm:main} implies that if $\lambda(G)=\sigma(G)=k$, then $G$ must be obtained by an iterated $k$-join of $k$-irreducible graphs.

\begin{proof} 
We start with $G$ and build the rooted tree $R$ and trees $T_H$ recursively.

Suppose that a graph $\Gamma$ is the label of a node which has yet to be considered.  If $\Gamma$ is $k$-irreducible, this node will be a leaf in $R$, and we add $\Gamma$ to $\mathcal{A}$.

If $\Gamma$ is $k$-reducible, then $\lambda(\Gamma)=\sigma(\Gamma)=k$, and so $\Gamma$ contains some collection of $k$-edge cuts.  By Lemma \ref{lemma:disjointcocircuits}, these must be pairwise non-overlapping.  Removing the edges from all of these $k$-edge cuts yields a graph with some number $d\geq 2$ of connected components; label these $\Gamma_1,\ldots,\Gamma_d$.  We observe that each $\Gamma_i$ must be $k$-edge connected and contains at least $k$ edge-disjoint spanning trees.

We note that $\Gamma$ is therefore an order-independent $k$-join of $\Gamma_1,\ldots,\Gamma_d$, and we can build a tree to specify explicitly which pairs are joined, as follows.  Define $T_\Gamma$ to be the graph obtained from $\Gamma$ by contracting each $\Gamma_i$ to a single vertex (and removing any multiple edges).  Lemma \ref{lemma:disjointcocircuits} ensures that this graph is a tree.  Now, each edge $e$ of $T_\Gamma$ corresponds to exactly one of the $k$-edge cuts of $\Gamma$, so we can label these $k$-edge cuts by the edges of $T_\Gamma$.

Finally, we add a child node of $\Gamma$ to $R$ labelled by $\Gamma_i$ for each $i \in \{1,\ldots,d\}$, and apply the recursion to each of these new nodes. 
\end{proof}

\begin{defn} \label{defn:ingredients}
For a given max-STP graph $G$, the {\em max-STP decomposition} of $G$ is the triple $\mathcal{I}(G) = (\mathcal{A},R,\{T_H\})$.
\end{defn}

\begin{example} \label{example:maxSTP}
Consider the max-STP graph $G$ shown in Figure \ref{figure:maxSTP}, which has $\lambda(G)=\sigma(G)=2$.  Now, $G$ possesses exactly one 2-edge cut, so the root vertex of $R$, labelled by $G$, has two descendants; the tree $T_G$ associated with it is the unique tree on two vertices.  Now, one of the child nodes is labelled by a copy of $K_4$, which is 2-irreducible, so the node becomes a leaf.  The other child node is labelled by a graph $H$ which is the order-independent 2-join of three copies of $K_4$.  Thus this node has three child nodes, all leaves labelled by a copy of $K_4$, and the associated tree $T_H$ is the unique tree on 3 vertices.

Thus the max-STP decomposition of $G$ are $(\mathcal{A},R,\mathcal{T} )$, where $\mathcal{A}$ contains of four copies of $K_4$, $R$ is as shown in Figure \ref{figure:maxSTP}, and \mbox{$\mathcal{T}= \{$ {\setlength{\unitlength}{7mm}
\begin{picture}(1.2,0.3)
\put(0.1,0.2){\blob}
\put(0.1,0.2){\line(1,0){1}}
\put(1.1,0.2){\blob}
\end{picture} },{\setlength{\unitlength}{7mm}
\begin{picture}(2.2,0.3)
\put(0.1,0.2){\blob}
\put(0.1,0.2){\line(1,0){1}}
\put(1.1,0.2){\blob}
\put(1.1,0.2){\line(1,0){1}}
\put(2.1,0.2){\blob}
\end{picture} } $\}$}.
\end{example}

\begin{figure}[hbtp]
\setlength{\unitlength}{7mm}
\centering
\begin{picture}(14,9)

\put(2.85,0){\complete}
\put(-0.15,4){\complete}
\put(2.85,7){\complete}
\put(5.85,4){\complete}

\put(3,2){\line(-1,2){1}}
\put(5,2){\line(1,2){1}}

\put(0,6){\line(1,1){3}}
\put(2,6){\line(1,1){1}}
\put(6,6){\line(-1,1){1}}
\put(8,6){\line(-1,1){3}}

\put(11,2){\blob}
\put(12,2){\blob}
\put(13,2){\blob}
\put(12,4){\blob}
\put(14,4){\blob}
\put(13,6){\blob}

\put(11,2){\line(1,2){1}}
\put(12,2){\line(0,1){2}}
\put(13,2){\line(-1,2){1}}
\put(12,4){\line(1,2){1}}
\put(14,4){\line(-1,2){1}}

\end{picture}
\caption{A max-STP graph $G$ with $\lambda(G)=\sigma(G)=2$, and the associated rooted tree $R$. \label{figure:maxSTP}}
\end{figure}
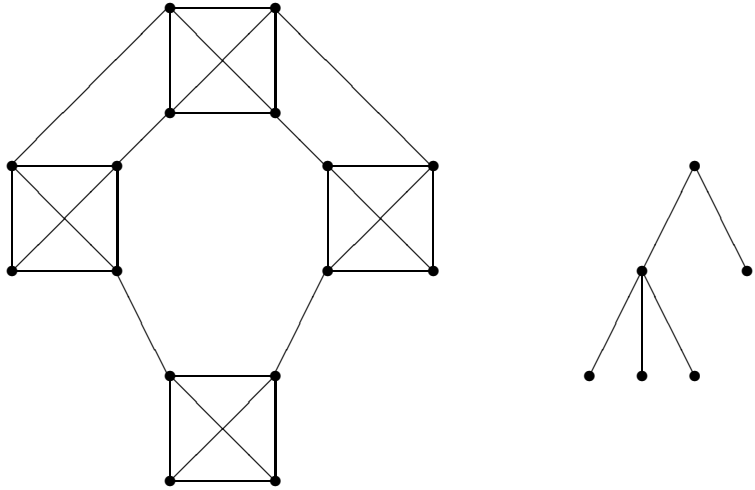

We remark that it is not possible to uniquely recover the original graph from a max-STP decomposition; the exact $k$-joins must be specified.  For example, the two graphs in Figure~\ref{figure:orderindep} are non-isomorphic, yet their max-STP decompositions are the same.  Moreover in any reconstruction of a graph from a max-STP decomposition, the collection of $k$-edge cuts of the graph at any node of $R$ must be precisely the set of $k$-joins applied to the graphs at the child nodes.  For example, none of the graphs in Figure~\ref{figure:possibilities} are valid reconstructions using the elements
\mbox{$( \{K_4,K_4,K_4,K_4\},\, R,\, \{$
{\setlength{\unitlength}{7mm}
\begin{picture}(1.2,0.3)
\put(0.1,0.2){\blob}
\put(0.1,0.2){\line(1,0){1}}
\put(1.1,0.2){\blob}
\end{picture} },{\setlength{\unitlength}{7mm}
\begin{picture}(2.2,0.3)
\put(0.1,0.2){\blob}
\put(0.1,0.2){\line(1,0){1}}
\put(1.1,0.2){\blob}
\put(1.1,0.2){\line(1,0){1}}
\put(2.1,0.2){\blob}
\end{picture} } $\})$}
of the max-STP decomposition from Example~\ref{example:maxSTP}; their actual max-STP decompositions are also given in Figure~\ref{figure:possibilities}.

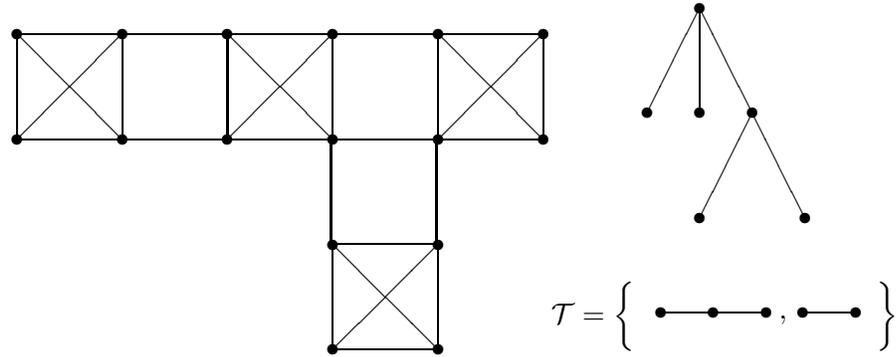
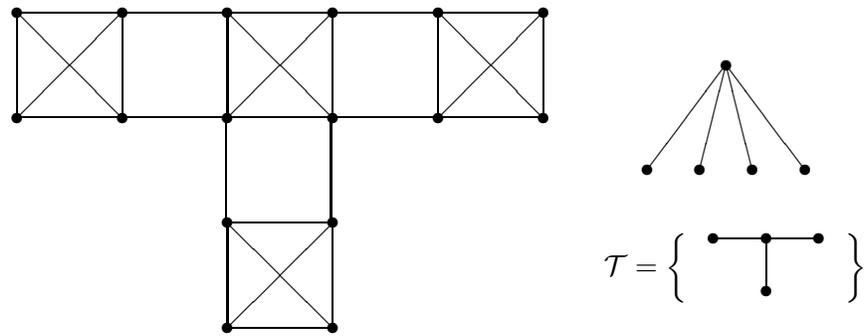
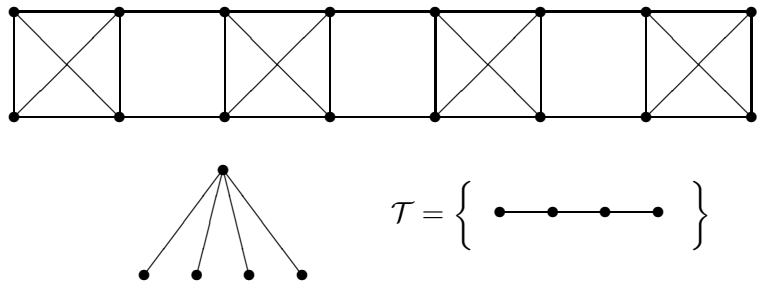
\begin{figure}[hbtp]
\setlength{\unitlength}{7mm}
\centering
\subfigure[Leaving some previous 2-joins order-independent of the new 2-join]{%
\begin{picture}(15,6.7)
\put(-0.15,4){\complete}
\put(3.85,4){\complete}
\put(7.85,4){\complete}
\put(5.85,0){\complete}

\put(2,4){\line(1,0){2}}
\put(2,6){\line(1,0){2}}
\put(6,4){\line(1,0){2}}
\put(6,6){\line(1,0){2}}
\put(6,2){\line(0,1){2}}
\put(8,2){\line(0,1){2}}

\put(12,4.5){\blob}
\put(13,2.5){\blob}
\put(13,4.5){\blob}
\put(13,6.5){\blob}
\put(14,4.5){\blob}
\put(15,2.5){\blob}

\put(12,4.5){\line(1,2){1}}
\put(13,4.5){\line(0,1){2}}
\put(14,4.5){\line(-1,2){1}}
\put(13,2.5){\line(1,2){1}}
\put(15,2.5){\line(-1,2){1}}

\put(10,0.5){ ${\displaystyle \mathcal{T}=\left\{ \phantom{\int} \right.}$
  \begin{picture}(3.0,0.3)
  \put(-0.3,0.2){\blob}
  \put(-0.3,0.2){\line(1,0){1}}
  \put(0.7,0.2){\blob}
  \put(0.7,0.2){\line(1,0){1}}
  \put(1.7,0.2){\blob}
  \put(1.95,0.1){,}
  \put(2.4,0.2){\blob}
  \put(2.4,0.2){\line(1,0){1}}
  \put(3.4,0.2){\blob}
  \end{picture} ${\displaystyle \left. \phantom{\int} \right\} }$ }

\end{picture} }

\subfigure[Leaving {\em all} previous 2-joins order-independent of the new 2-join]{%
\begin{picture}(15,6.7)
\put(-0.15,4){\complete}
\put(3.85,4){\complete}
\put(7.85,4){\complete}
\put(3.85,0){\complete}

\put(2,4){\line(1,0){2}}
\put(2,6){\line(1,0){2}}
\put(6,4){\line(1,0){2}}
\put(6,6){\line(1,0){2}}
\put(4,2){\line(0,1){2}}
\put(6,2){\line(0,1){2}}

\put(12,3){\blob}
\put(13,3){\blob}
\put(14,3){\blob}
\put(15,3){\blob}
\put(13.5,5){\blob}

\put(12,3){\line(3,4){1.5}}
\put(13,3){\line(1,4){0.5}}
\put(14,3){\line(-1,4){0.5}}
\put(15,3){\line(-3,4){1.5}}

\put(11,1){ ${\displaystyle \mathcal{T}=\left\{ \phantom{\int} \right.}$
  \begin{picture}(1.4,1)
  \put(-0.3,0.7){\blob}
  \put(0.7,0.7){\blob}
  \put(1.7,0.7){\blob}
  \put(0.7,-0.3){\blob}
  \put(-0.3,0.7){\line(1,0){1}}
  \put(0.7,0.7){\line(1,0){1}}
  \put(0.7,-0.3){\line(0,1){1}}
  \end{picture} ${\displaystyle \left. \phantom{\int} \right\} }$ }

\end{picture} }

\subfigure[Same as (b), but with a different set $\mathcal{T}$]{%
\begin{picture}(14,5.7)
\put(-0.15,3){\complete}
\put(3.85,3){\complete}
\put(7.85,3){\complete}
\put(11.85,3){\complete}

\put(2,3){\line(1,0){2}}
\put(2,5){\line(1,0){2}}
\put(6,3){\line(1,0){2}}
\put(6,5){\line(1,0){2}}
\put(10,3){\line(1,0){2}}
\put(10,5){\line(1,0){2}}

\put(2.5,0){\blob}
\put(3.5,0){\blob}
\put(4.5,0){\blob}
\put(5.5,0){\blob}
\put(4,2){\blob}

\put(2.5,0){\line(3,4){1.5}}
\put(3.5,0){\line(1,4){0.5}}
\put(4.5,0){\line(-1,4){0.5}}
\put(5.5,0){\line(-3,4){1.5}}

\put(7,1){ ${\displaystyle \mathcal{T}=\left\{ \phantom{\int} \right.}$
\begin{picture}(2.5,0.3)
  \put(-0.3,0.2){\blob}
  \put(-0.3,0.2){\line(1,0){1}}
  \put(0.7,0.2){\blob}
  \put(0.7,0.2){\line(1,0){1}}
  \put(1.7,0.2){\blob}
  \put(1.7,0.2){\line(1,0){1}}
  \put(2.7,0.2){\blob}
\end{picture} ${\displaystyle \left. \phantom{\int} \right\} }$ }

\end{picture} }

\caption{Three graphs constructed from the same max-STP decomposition, but each with different max-STP decompositions. 
\label{figure:possibilities}}
\end{figure}


\section{Decomposing matroids} \label{section:matroiddecomp}

From now on, all matroids we consider will be connected.  In the case of matroids with $k$ disjoint bases, we introduce the following terminology which is analogous to the notion of $k$-reducibility for graphs we saw above.
\newpage
\begin{defn} \label{defn:reducible}
Suppose that $M$ is a matroid with $\sigma(M) \geq k$.  We call $M$ {\em $k$-reducible} if $\lambda(M)=\sigma(M)=k$; otherwise (i.e.~if either $\sigma(M)>k$ or $k=\sigma(M)<\lambda(M)$) we call $M$ {\em $k$-irreducible}.
\end{defn}

We note that this partitions the class of matroids with $k$ disjoint bases (as characterized by Edmonds \cite{Edmonds65}) into either being $k$-reducible or $k$-irreducible.  We also remark that the case of isolated vertices does not arise here, as their cycle matroids are empty.

As the object obtained from a matroid $M$ by deleting all of its minimum cocircuits will be keep appearing, it is useful to give it a formal name.

\begin{defn} \label{defn:crux}
Let $M$ be a matroid whose minimum cocircuits are $C_1,\ldots,C_\ell$.  Then the {\em crux} of $M$, denoted $\chi(M)$, is defined to be
\[ \chi(M) = M\setminus \bigcup_{i=1}^\ell C_i. \]
\end{defn}

If it happens that $\bigcup_{i=1}^\ell C_i = E$, then by abuse of notation we write $\chi(M)=\emptyset$ (as a shorthand for the matroid $(\emptyset,\{\emptyset\})$).

\begin{lemma} \label{lemma:cocircuits}
Suppose that $\lambda(M)=\sigma(M)=k$, and that $\{C_1,\ldots,C_\ell\}$ are the minimum cocircuits of $M$.  Then for all $i \neq j$, $C_i$ is a minimum cocircuit of $M\setminus C_j$.
\end{lemma}

\proof Suppose that $C_i$ and $C_j$ are minimum cocircuits of $M$, and that $\rank(M)=r$.  Their complements are hyperplanes, so the intersection of their complements is a flat of rank at most $r-2$ (since it is strictly contained in two distinct hyperplanes).  Since $\lambda(M)=\sigma(M)$, we have a base that intersects $C_i$ and $C_j$ in exactly one element each, so by removing these elements we have an independent set of size $r-2$ in the intersection of the two hyperplanes.  Therefore, $M\setminus C_j$ is a matroid of rank $r-1$, which contains a flat $A$ of rank $r-2$, whose complement is $C_i$.  Consequently, $C_i$ is a cocircuit of $M\setminus C_j$.  

To show that $C_i$ is minimum, it suffices to show that $\lambda(M\setminus C_j)\geq k$.  Since $M$ has $k$ disjoint bases, each of which intersects $C_j$ in a single element, we have $\sigma(M\setminus C_j)\geq k$, and thus $\lambda(M\setminus C_j)\geq k$.  Hence $C_i$ is a minimum cocircuit of $M\setminus C_j$.
\endproof

\begin{lemma} \label{lemma:cruxrank}
Suppose that $M$ is a matroid such that $\lambda(M)=\sigma(M)=k$, $\rank(M)=r$, and whose minimum cocircuits are $C_1,\ldots,C_\ell$.  Then the crux of $M$ has rank $r-\ell$.
\end{lemma}

\proof By applying Lemma~\ref{lemma:cocircuits} repeatedly, we see that deleting each minimum cocircuit in turn will reduce the rank by~$1$.  Once we have deleted all minimum cocircuits $C_1,\ldots,C_\ell$ to obtain the crux, we have a matroid of rank $r-\ell$.  \endproof

Even if a matroid $M$ is assumed to be connected, its crux $\chi(M)$ is not necessarily connected.  (For instance, in the case of graphic matroids, this is clear.)  So the following parameter makes sense.

\begin{notation}
We let $\delta(M)$ denote the number of connected components of $\chi(M)$.
\end{notation}

Note that if $\chi(M)=\emptyset$, we define $\delta(M)=0$. We remark that if $\sigma(M)=\lambda(M)=k$ and $\chi(M)=\emptyset$, then the matroid $M$ is the disjoint union of the minimum cocircuits, the bases are precisely transversals of the partition into minimum cocircuits, and thus $M$ is disconnected.  (These are precisely the transversal matroids discussed at the end of Section~\ref{section:intro}.)  Hence, if we restrict ourselves to connected matroids, the possibility that the crux is empty does not arise.

The next definition allows us to encode how a matroid with $\sigma(M)=\lambda(M)=k$ is assembled from its cocircuits and the connected components of its crux.

\begin{defn} \label{defn:assembly}
Let $M$ be a connected matroid such that $\sigma(M)=\lambda(M)=k$, with~$\ell$ minimum cocircuits $C_1,\ldots,C_\ell$, and where $\chi(M)$ has $d$ connected components $K_1,\ldots,K_d$.  For each minimum cocircuit $C_j$, let $\mathcal{V}_j$ denote the largest subset of $\{K_1,\ldots,K_d\}$ such that the restriction of $M$ to 
\[ C_j \cup \left( \bigcup_{K\in\mathcal{V}_j} K \right) \]
is connected.
Then the {\em assembly hypergraph} of $M$, denoted $\mathcal{H}(M)$, is the non-uniform hypergraph whose vertices are labelled by the connected components $K_1,\ldots,K_d$, where the hyperedges are labelled by the cocircuits $C_1,\ldots,C_\ell$, and the vertices incident with $C_j$ are precisely the members of $\mathcal{V}_j$.
\end{defn}

In other words, the assembly hypergraph $\mathcal{H}(M)$ tell us which components of $\chi(M)$ are ``joined'' by each cocircuit $C_j$ in the matroid $M$.

\begin{thm} \label{thm:decomp}
Suppose that $M=(E,\mathcal{I})$ is a matroid for which $\sigma(M)=\lambda(M)=k$.  Then we have the following.
\begin{enumerate}
\item There exists a unique set of $k$-irreducible matroids $\mathcal{M}=\{M_1,\ldots,M_m\}$ (for some integer $m$).
\item There exists a unique rooted tree $R$ with $m$ leaves labelled by $M_1,\ldots,M_m$, such that the root is labelled by $M$ and each non-leaf labelled by $K$ has $d=\delta(K)$ children, labelled by the connected components of $\chi(K)$.
\item For each non-leaf, labelled by $K$ and its $d$ children labelled $K_1,\ldots,K_d$, there exists a unique assembly hypergraph with $\ell$ hyperedges, and where $\sum_{i=1}^d \rank(K_i) = \rank(K)-\ell$.
\end{enumerate}
\end{thm}

\proof We build the rooted tree $R$ recursively, obtaining $\mathcal{M}$ and the collection of assembly matrices as we go along.  We begin by assigning the root node to our matroid $M$, and declaring $\mathcal{M}=\emptyset$.

Suppose that a matroid $K$ is the label of a node which is yet to be considered.  If $K$ is $k$-irreducible, then this node becomes a leaf in $R$, and we add $K$ to $\mathcal{M}$.

On the other hand, if $K$ is $k$-reducible, then $\lambda(K)=\sigma(K)=k$, and so $K$ has some minimum cocircuits of size $k$; suppose that there are $\ell$ of these, labelled $C_1,\ldots,C_\ell$.  By Lemma \ref{lemma:disjointcocircuits}, these are all disjoint.  Now consider the crux $\chi(K)$; by repeatedly applying Lemma \ref{lemma:deletecocircuit}, this must have at least $k$ disjoint bases.  Now suppose that the crux has $d=\delta(K)$ connected components $K_1,\ldots,K_d$; since $K$ is connected, there must be at least one of these.  Furthermore, since a base for $\chi(K)$ is the disjoint union of a base for each component, each $K_i$ has $\sigma(K_i)\geq k$, and thus each $K_i$ must have a cocircuit of size at least $k$.  By Lemma \ref{lemma:cruxrank}, the rank of $\chi(K)$ is $\rank(K)-\ell$, so the ranks of its connected components $K_1,\ldots,K_d$ must sum to this.

For each $j\in\{1,\ldots,\ell\}$, we add a hyperedge to the assembly hypergraph $\mathcal{H}(M)$ as follows.  Consider the matroid $L$ formed by restricting $K$ to $\chi(K)\cup C_j$.  Now, the component of $L$ containing $C_j$ will be of the form 
\[ C_j \cup \left( \bigcup_{i \in I} K_i \right), \]
where $I \subseteq \{1,\ldots,d\}$.  Then we add the hyperedge $\{ K_i \mid i \in I\}$.

Finally, we add a child node of $K$ to $R$ labelled by $K_i$ for each $i \in \{1,\ldots,d\}$, and apply the recursion to each of these new nodes.  \endproof

\section{Discussion} \label{section:discussion}

We conclude the paper with a few remarks about our main results (Theorem~\ref{thm:main} for graphs and Theorem~\ref{thm:decomp} for matroids).  We also consider some issues related to computational complexity.

A general comment about Theorem~\ref{thm:main} is perhaps in order.  The graph clearly determines the
decomposition, but the converse is not true.  In general, it is not possible to recover the
graph from the decomposition, except in some special cases (for instance, if the original graph
is a tree).  In general, even for $k=1$ we cannot recover the graph: if the $k$-irreducible
subgraphs have more than one vertex there are many different ways the edge-cuts could be added
back in.  While of course one could keep additional information at each stage in order to
permit this reconstruction, we prefer to view these results as a general structural description
of the original graph rather than as a precise encoding of it.  In particular, Theorem~\ref{thm:main}
shows that, at each stage of the decomposition, a max-STP graph $G$ is globally ``tree-like'', in that 
it may be contracted to the tree $T_G$.  A similar comment applies to Theorem~\ref{thm:decomp}.

We also observe that there are some constraints on the ingredients of a max-STP decomposition.
For the case of graphs (Theorem~\ref{thm:main}), consider for definiteness the root node of $R$ and
one of its child nodes.  The root is labelled by $G$ and has an associated tree $T_G$.  The
child is labelled by $H$ and has an associated tree $T_H$.  One of the vertices of $T_G$, say $v_H$,
corresponds to $H$; let the degree of $v_H$ be $d$.  Then there are $d$ $k$-edge cuts in
$G$ incident with $H$.  This collection of edge-cuts must be sufficient so that no $k$-edge cut of 
$H$ is an edge-cut of $G$.  
One way of characterizing this condition is that for every edge $e$ of $T_H$ there must exist a pair of edges 
from one of the $d$ $k$-joins connected to $H$ in $T_G$ such that each edge is incident to a different
connected component of $T_H \setminus e$.
This gives a constraint on the structure of $T_H$.  For instance it is
not hard to see from these constraints that if $k=1$ then $R$ has a root node labelled by
$T_G$ and child nodes labelled by the $k$-irreducible subgraphs.  Note also that if a node
of $R$ has an associated tree then the number of children of that node is the number of
vertices of the associated tree; otherwise that node of $R$ is a leaf labelled by a
$k$-irreducible graph.

It is important to note that Theorem~\ref{thm:main} is not a simple corollary of Theorem~\ref{thm:decomp}.
Given a graph $G$, we may apply Theorem~\ref{thm:main} directly to $G$ or we may apply
Theorem~\ref{thm:decomp} to the matroid $M(G)$; however, we get two different decompositions.  In the
latter the ``connected'' components are now blocks of the underlying graph, not connected
components in the graph-theoretic sense.

Finally, we remark that $\sigma$ can be computed in polynomial time for both graphs and
matroids (see Schrijver \cite[Sections~51.4 and~42.3]{Schrijver03}).  As $\lambda$ can also be
found in polynomial time for graphs (\cite[Section~15.3]{Schrijver03}) this means that we can
determine if a graph is max-STP in polynomial time.  Furthermore, it is actually possible to
find the collection of $\sigma(G)$ edge-disjoint spanning trees in polynomial time; the best
algorithm known for doing this is due to Gabow and Westermann \cite{Gabow92}.  However for
matroids in general, or even binary matroids, determining the girth (or cogirth) was shown to
be NP-hard by Cho {\em et al.}\ in \cite{Cho07}: this follows from the equivalence of the
problem to that of determining the minimum distance of a binary linear code, which was shown to
be NP-hard by Vardy \cite{Vardy97}.  Thus, unless $\rm{P}=\rm{NP}$ there is no obvious polynomial time
algorithm for testing $\lambda(M)=\sigma(M)$ for an arbitrary matroid $M$.

\section*{Acknowledgements}
The authors acknowledge financial support from NSERC and the Ontario Ministry of Research and Innovation.

\end{document}